\documentclass[12pt,a4paper,thmsb]{article}
\usepackage{sw20lart}

\input tcilatex
\newcommand{\R}{{\sf R\hspace*{-0.9ex}\rule{0.15ex}{1.5ex}\hspace*{0.9ex}}}
\begin{document}

\title{K-Bessel functions associated \\
to 3-rank Jordan algebra.}
\author{\textsc{Hacen DIB} \\
Department of Mathematics, University of Tlemcen\\
B.P. 119 Tlemcen 13000, ALGERIA\\
e-mail : h\_dib@mail.com}
\date{}
\maketitle

\begin{abstract}
Using Bessel-Muirhead system, we can express the K-bessel function defined
on a Jordan algebra as linear combination of the J-solutions. We determine
explicitly the coefficients when the rank of this Jordan algebra is three
after a reduction to the rank two. The main tools are some algebraic
identities developed for the occasion.
\end{abstract}



\section{Introduction}

In \cite{dib1} we solved the Bessel-Muirhead system of rank 2 and 3 and
expressed, in the rank 2 case, the K-bessel function as linear combination
of the J-solutions with explicit coefficients. In this paper we continue our
work and prove that similar result is also true when the Jordan algebra is
of rank 3. In this case there is four non-equivalent classes of simple
euclidean Jordan algebra and in \cite{dib1} we intended to perform
case-by-case calculation. In this way, a serious difficulty arise in the
evaluation of some integral over the automorphism group of the Jordan
algebra. However, a unified treatment is possible by a reduction to the rank
2 case using some algebraic identities. This what we will present here. Let
us recall the situation and some results in \cite{dib1}.

\begin{definition}
Bessel-Muirhead operators are defined by 
\begin{equation}
B_i=x_i\frac{\partial ^2}{\partial x_i^2}\;+\;(\nu +1)\frac \partial
{\partial x_i}\;+\;1\;+\;\frac d2\stackunder{j\neq i}{\sum }\frac
1{x_i-x_j}(x_i\frac \partial {\partial x_i}\;-\;x_j\frac \partial {\partial
x_j})\;,\;1\leq i\leq r  \label{muirh}
\end{equation}
$r$ \ is the rank of the system. A symmetric function $\;f$ \ is said to be
a Bessel function if it is a solution of \ $B_i\;f=0\;,i=1,2,...,r.$
\end{definition}

Denote by $\;t_1,t_2,...,t_r$ \ the elementary symmetric functions , i.e.

\begin{equation}
t_p=\stackunder{1\leq i_1\leq i_2\leq ...\leq i_p\leq r}{\dsum }%
\;x_{i_1}x_{i_2}...x_{i_p}
\end{equation}
with $t_0=1$ and $t_p=0$ \ if $\;p<0$ \ or \ $p>r$. The Bessel-Muirhead
system is then equivalent to the system (see \cite{james},\cite{dib2}) $%
\;Z_kg=0\;,\;1\leq k\leq r$ \ where

\begin{equation}
Z_k=\stackunder{i,j=1}{\stackrel{r}{\sum }}A_{ij}^k\frac{\partial ^2}{%
\partial t_i\partial t_j}+(\nu +1+\frac{r-k}2d)\frac \partial {\partial
t_k}+\delta _k^1  \label{mod}
\end{equation}
and 
\begin{equation}
A_{ij}^k=\left\{ 
\begin{array}{ccc}
t_{i+j-k} & if & i,j\geq k \\ 
-t_{i+j-k} & if & i,j<k\;and\;i+j\geq k \\ 
0 &  & elswhere
\end{array}
\right.
\end{equation}
Here \ $\delta _k^1\;$is the Kronecker symbol and $%
g(t_1,t_2,...,t_r)=f(x_1,x_2,...,x_r).$ When $r=2$, we have (see \cite{dib1}
) a fundamental set of solutions given by $J_\nu ^{\left[ 2,1\right] },J_\nu
^{\left[ 2,2\right] },t_2^{-\nu }J_{-\nu }^{\left[ 2,1\right] }$ and $%
t_2^{-\nu }J_{-\nu }^{\left[ 2,2\right] }$ where : 
\begin{equation}
J_\nu ^{\left[ 2,1\right] }(t_1,t_2)=\stackunder{m_1,m_2\geq 0}{\dsum }\frac{%
(-1)^{m_1}}{(1)_{m_1}(1)_{m_2}(1+\nu )_{m_2}(1+\nu +\frac d2)_{m_1+2m_2}}%
t_1^{m_1}t_2^{m_2}
\end{equation}
and 
\begin{equation}
J_\nu ^{\left[ 2,2\right] }(t_1,t_2)=\stackunder{m_1,m_2\geq 0}{\dsum }\frac{%
(-1)^{m_1}}{(1-\nu -d/2)_{m_1}(1)_{m_2}(1+\nu )_{m_2}(1)_{m_1+2m_2}}%
t_1^{m_1-\nu -d/2}t_2^{m_2}
\end{equation}
Here $(a)_k$ is the classical Pochhammer symbol i.e, $(a)_k=a(a+1)\cdots
(a+k-1).$ In this case the K-Bessel function (in two variables) can be
written as follows : 
\begin{eqnarray}
K_\nu ^{\left[ 2\right] } &=&\left( 2\pi \right) ^{(n-2)/2}\Gamma (-\nu
)\Gamma (-\nu -\frac{n-2}2)J_\nu ^{\left[ 2,1\right] }+\left( 2\pi \right)
^{(n-2)/2}\Gamma (\nu )\Gamma (\nu -\frac{n-2}2)t_2^{-\nu }J_{-\nu }^{\left[
2,1\right] } \\
&&\ +\left( 2\pi \right) ^{(n-2)/2}\Gamma (-\nu )\Gamma (\nu +\frac{n-2}%
2)J_\nu ^{\left[ 2,2\right] }+\left( 2\pi \right) ^{(n-2)/2}\Gamma (\nu
)\Gamma (-\nu +\frac{n-2}2)t_2^{-\nu }J_{-\nu }^{\left[ 2,2\right] } 
\nonumber
\end{eqnarray}
Observe that this expression reduce (up to constant factor) to the classical
one variable formula when $t_2=0$ and generic $\nu $. Now, for $r=3$ we have
eight linearly independant J-solutions $J_\nu ^{\left[ 3,1\right] },J_\nu
^{\left[ 3,2\right] },J_\nu ^{\left[ 3,3\right] },J_\nu ^{\left[ 3,4\right]
},t_3^{-\nu }J_{-\nu }^{\left[ 3,1\right] },t_3^{-\nu }J_{-\nu }^{\left[
3,2\right] },t_3^{-\nu }J_{-\nu }^{\left[ 3,3\right] }$ and $t_3^{-\nu
}J_{-\nu }^{\left[ 3,4\right] }$ where : 
\begin{eqnarray}
J_\nu ^{\left[ 3,1\right] }(t_1,t_2,t_3) &=&\stackunder{m_1,m_2,m_3\geq 0}{%
\dsum }\frac{(-1)^{m_1+m_3}\;}{(1)_{m_1}(1)_{m_2}(1)_{m_3}(1+\nu
)_{m_3}(1+\nu +\frac d2)_{m_2+2m_3}}\times \\
&&\times \frac{(1+2\nu +d)_{m_1+2m_2+4m_3}}{(1+\nu
+d)_{m_1+2m_2+3m_3}(1+2\nu +d)_{m_1+2m_2+3m_3}}\;t_1^{m_1}t_2^{m_2}t_3^{m_3}
\nonumber
\end{eqnarray}
\begin{eqnarray}
J_\nu ^{\left[ 3,2\right] }(t_1,t_2,t_3) &=&\stackunder{m_1,m_2,m_3\geq 0}{%
\dsum }\frac{(-1)^{m_1+m_3}\;}{(1-\nu -d)_{m_1}(1)_{m_2}(1)_{m_3}(1+\nu
)_{m_3}(1+\nu +\frac d2)_{m_2+2m_3}}\;\times \\
&&\times \frac{(1+\nu )_{m_1+2m_2+4m_3}}{(1)_{m_1+2m_2+3m_3}(1+\nu
)_{m_1+2m_2+3m_3}}\;t_1^{m_1-\nu -d}t_2^{m_2}t_3^{m_3}  \nonumber
\end{eqnarray}
\begin{eqnarray}
J_\nu ^{\left[ 3,3\right] }(t_1,t_2,t_3) &=&\stackunder{m_1,m_2,m_3\geq 0}{%
\dsum }\frac{(-1)^{m_1+m_3}}{(1)_{m_1}(1-\nu -d/2)_{m_2}(1)_{m_3}(1+\nu
)_{m_3}(1)_{m_2+2m_3}}\;\times \\
&&\times \frac{(1)_{m_1+2m_2+4m_3}}{(1-\nu
)_{m_1+2m_2+3m_3}(1)_{m_1+2m_2+3m_3}}\;t_1^{m_1}t_2^{m_2-\nu -d/2}t_3^{m_3} 
\nonumber
\end{eqnarray}
\begin{eqnarray}
J_\nu ^{\left[ 3,4\right] }(t_1,t_2,t_3) &=&\stackunder{m_1,m_2,m_3\geq 0}{%
\dsum }\frac{(-1)^{m_1+m_3}}{(1+\nu )_{m_1}(1-\nu -d/2)_{m_2}(1)_{m_3}(1+\nu
)_{m_3}(1)_{m_2+2m_3}}\;\times \\
&&\times \frac{(1+\nu )_{m_1+2m_2+4m_3}}{(1)_{m_1+2m_2+3m_3}(1+\nu
)_{m_1+2m_2+3m_3}}\;t_1^{m_1+\nu }t_2^{m_2-\nu -d/2}t_3^{m_3}  \nonumber
\end{eqnarray}
Observe also that when $t_3=0$ (and $\nu $ generic) these functions reduce
to : 
\begin{equation}
J_\nu ^{\left[ 3,1\right] }(t_1,t_2,0)=J_{\nu +d/2}^{\left[ 2,1\right]
}(t_1,t_2)  \label{red1}
\end{equation}
\begin{equation}
J_\nu ^{\left[ 3,2\right] }(t_1,t_2,0)=J_{\nu +d/2}^{\left[ 2,2\right]
}(t_1,t_2)  \label{red2}
\end{equation}
\begin{equation}
J_\nu ^{\left[ 3,3\right] }(t_1,t_2,0)=t_2^{-\nu -d/2}\;J_{-\nu
-d/2}^{\left[ 2,1\right] }(t_1,t_2)  \label{red3}
\end{equation}
\begin{equation}
J_\nu ^{\left[ 3,4\right] }(t_1,t_2,0)=t_2^{-\nu -d/2}\;J_{-\nu
-d/2}^{\left[ 2,2\right] }(t_1,t_2)  \label{red4}
\end{equation}

\section{Some algebraic identities}

For the general theory of Jordan algebra one can see \cite{fk}, but what we
will develop is somehow specific to the rank three. So, let $A$ a real
simple and euclidean Jordan algebra with rank $3$ and dimension $n$. We know
tha $n=3+3d$ where $d=1,2,4$ or $8$. Let $\left\{ c_1,c_2,c_3\right\} $ be a
complete system of mutually orthogonal primitive idempotents i.e, $%
c_ic_j=\delta _i^jc_i$, $c_1+c_2+c_3=e$ the unit of $A$ and none of the $c_j$%
's can split into a sum of two idempotents. We have a Cayley-Hamilton like
theorem: $x^3-a_1(x)x^2+a_2(x)x-a_3(x)e=0$ and a spectral decomposition: $%
x=k.\left( \lambda _1c_1+\lambda _2c_2+\lambda _3c_3\right) $ with $k$ an
element of the automorphism group of $A$ and $\lambda _i$ reals such that: 
\begin{eqnarray}
a_1(x) &=&\lambda _1+\lambda _2+\lambda _3:=tr(x)  \nonumber \\
a_3(x) &=&\lambda _1\lambda _2\lambda _3:=\det (x)\text{\ and} \\
a_2(x) &=&\lambda _1\lambda _2+\lambda _1\lambda _3+\lambda _2\lambda
_3=\frac 12\left[ tr(x)^2-tr(x^2)\right]   \nonumber
\end{eqnarray}
The inner product is defined then by: $\left( x,y\right) :=tr(xy)$. The
operators $L(x)$ and $P(x)$ are defined by $L(x)y=xy$ and $%
P(x)=2L^2(x)-L(x^2)$. Let us consider the Peirce decomposition with respect
to the idempotent $c_3$ i.e, $A=A_0\oplus A_{1/2}\oplus A_1$ where $A_\alpha 
$ is the eigenspace of $L(c_3)$ with respect to the eigenvalue $\alpha $. $%
A_0$ and $A_1$ are Jordan subalgebras of rank 2 and 1 respectively, and $%
A_{1/2}$ is a subspace of dimension $2d$. Put $n_0=\dim A_0=2+d$ and $%
e_0=c_1+c_2$ the unit of $A_0$. We have $A_0A_1=\left\{ 0\right\} $, $\left(
A_0\oplus A_1\right) A_{1/2}\subset A_{1/2}$ and $A_{1/2}A_{1/2}\subset
A_0\oplus A_1$. If we write $tr(z)$ or $\det (z)$ of an element of $A_0$
this will mean trace and determinant with respect to the subalgebra $A_0$.
We denote by $\Omega _3$ the cone of positivity of $A$ i.e, $\Omega
_3=\left\{ x\in A\;/\;\lambda _i>0,i=1,2,3\right\} =\left\{ x\in
A\;/\;a_i(x)>0,i=1,2,3\right\} $ and by $\Omega _2$ the cone of $A_0$. Every 
$x$ in $\Omega _3$ (resp. in $\Omega _2$) admit a unique square root in $%
\Omega _3$ (resp. in $\Omega _2$) and is invertible.

\begin{lemma}
For $y=e_0+\xi +tc_3$ with $\xi \in A_{1/2}$ and $t\in \TeXButton{R}{{\R}}$
we have 
\begin{equation}
\det (y)=t-\frac 12\left\| \xi \right\| ^2  \label{lem1}
\end{equation}
\end{lemma}

\textbf{proof} : The projection onto $A_1$ is $P\left( c_3\right) $, so $%
P\left( c_3\right) \xi =0$ and therefore $0=tr\left( P\left( c_3\right) \xi
\right) =\left( e,P\left( c_3\right) \xi \right) =\left( c_3,\xi \right)
=\frac 12tr(\xi )$. By the same argument $tr(\xi ^3)=0$. Now by
Cayley-Hamilton $\xi ^3+a_2(\xi )\xi -\det (\xi )e=0$ which implies $\det
(\xi )=0$. On the other hand $\xi ^2=u+\tau c_3$ with $u\in A_0$ and $\tau
\in \TeXButton{R}{{\R}}$. We have $\tau =\left( c_3,\xi ^2\right) =\left(
\xi c_3,\xi \right) =\frac 12\left\| \xi \right\| ^2$ and $tr(u)=\frac
12\left\| \xi \right\| ^2$. From $\det (\xi ^2)=\tau \det (u)$ we deduce
that $\det (u)=0$. So by Cayley-Hamilton (in $A_0$) we can write $u^2=\frac
12\left\| \xi \right\| ^2u$. Now $\xi ^2=u+\tau c_3\Rightarrow \xi
^4=u^2+\tau ^2c_3=\frac 12\left\| \xi \right\| ^2u+\frac 14\left\| \xi
\right\| ^4c_3=\frac 12\left\| \xi \right\| ^2\xi ^2$ and then $\xi ^3=\frac
12\left\| \xi \right\| ^2\xi $. Therefore $u\xi =\frac 14\left\| \xi
\right\| ^2\xi $. Then we have : 
\[
y=e_0+\xi +tc_3 
\]
\[
y^2=\left( e_0+u\right) +\left( 1+t\right) \xi +\left( t^2+\frac 12\left\|
\xi \right\| ^2\right) c_3 
\]
\[
y^3=\left( e_0+\left( 2+t\right) u\right) +\left( 1+t+t^2+\frac 12\left\|
\xi \right\| ^2\right) \xi +\left( t^2+t\left\| \xi \right\| ^2+\frac
12\left\| \xi \right\| ^2\right) c_3 
\]
and also $tr(y)=2+t$, $a_2(y)=\dfrac 12\left[ tr(y)^2-tr(y^2)\right]
=1+2t-\dfrac 12\left\| \xi \right\| ^2$. The result is a consequence of 
\[
y^3-tr(y)y^2+a_2(y)y=\left( t-\frac 12\left\| \xi \right\| ^2\right) e 
\]

\begin{lemma}
If $y=z+\xi +tc_3$ with $z\in \Omega _2,\;\xi \in A_{1/2}$ and $t\in 
\TeXButton{R}{{\R}}$ we have 
\begin{equation}
\det (y)=\det (z)\left[ t-\left( z^{-1}\xi ,\xi \right) \right]  \label{lem2}
\end{equation}
\end{lemma}

\textbf{proof} : The application 
\begin{eqnarray}
\rho &:&A_0\longrightarrow End\left( A_{1/2}\right) \\
u &\longrightarrow &\rho (u)\text{ defined by }\rho (u)\xi =2L(u)\xi 
\nonumber
\end{eqnarray}
is a representation of $A_0$ in the space $A_{1/2}$ (for more details see 
\cite{fk}). This mean that \\$2\rho (uv)=\rho (u)\rho (v)+\rho (v)\rho (u)$
and identically $P(u)\xi =0$. So 
\begin{eqnarray}
P(u+c_3)y &=&P(u+c_3)z+P(u+c_3)\xi +tP(u+c_3)c_3  \nonumber \\
&=&P(u)z+2L(u)\xi +tc_3
\end{eqnarray}
Now if $u=z^{-1/2}$ we derive the desired result thanks to $\left( \ref{lem1}%
\right) $ and the fact that 
\begin{eqnarray}
\det \left[ P(z^{-1/2}+c_3)y\right] &=&\det (z^{-1/2}+c_3)^2\det (y) 
\nonumber \\
&=&\det (z^{-1})\det (y) \\
&=&\det (z)^{-1}\det (y)  \nonumber
\end{eqnarray}

\begin{corollary}
If $y=z+\xi +tc_3$ $\in \Omega _3$, then : 
\begin{equation}
tr\left( y^{-1}\right) =\frac{2\det (z)+2ttr(z)-\left\| \xi \right\| ^2}{%
2\det (z)\left[ t-\left( z^{-1}\xi ,\xi \right) \right] }
\end{equation}
\end{corollary}

\textbf{proof} : first we have $y^3-tr(y)y^2+a_2(y)y-\det (y)e=0$. Then 
\[
y^{-1}=\dfrac 1{\det (y)}\left[ y^2-tr(y)y+a_2(y)e\right] 
\]
and therefore 
\begin{eqnarray*}
tr\left( y^{-1}\right) &=&\dfrac 1{\det (y)}\left[ tr\left( y^2\right)
-tr(y)^2+3a_2(y)\right] \\
&=&\dfrac 1{2\det (y)}\left[ tr(y)^2-tr\left( y^2\right) \right] \\
&=&\dfrac 1{2\det (y)}\left[ tr(z)^2+2ttr(z)-tr(z^2)-\left\| \xi \right\|
^2\right] \\
&=&\dfrac 1{2\det (y)}\left[ 2\det (z)+2ttr(z)-\left\| \xi \right\| ^2\right]
\end{eqnarray*}

\section{K-Bessel function}

The K-Bessel function is defined by (see \cite{dib2},\cite{fk})

\begin{equation}
K_\nu ^{\left[ 3\right] }(x)=\stackunder{\Omega _3}{\int }%
e^{-tr(y^{-1})-\left( x,y\right) }(\det y)^{\nu -\frac n3}dy
\end{equation}
After a change of variable, one can show that

\begin{equation}
K_\nu ^{\left[ 3\right] }(x)=(\det x)^{-\nu }K_{-\nu }^{\left[ 3\right] }(x)
\label{sym}
\end{equation}
Following \cite{dib2} where it is proved that $K_\nu ^{\left[ 3\right] }$ is
a solution of a differential system similar to (\ref{muirh}), we can write

\begin{equation}
K_\nu ^{\left[ 3\right] }(x)=\stackunder{j=1}{\stackrel{4}{\dsum }}a_\nu
^jJ_\nu ^{\left[ 3,j\right] }+b_\nu ^jt_3^{-\nu }J_{-\nu }^{\left[
3,j\right] }
\end{equation}
According to (\ref{sym}) we have : $a_\nu ^j=b_{-\nu }^j$ for $j=1,2,3,4$.
For suitable $\nu $, the following limit holds (see \cite{fk} for more
information on $\Gamma _{\Omega _3}$, the gamma function of the cone $\Omega
_3$) :

\begin{equation}
\stackunder{
\begin{array}{c}
x\rightarrow 0 \\ 
x\in \Omega _3
\end{array}
}{\lim }K_\nu ^{\left[ 3\right] }(x)=\Gamma _{\Omega _3}(-\nu )=(2\pi )^{%
\frac{3d}2}\Gamma (-\nu )\Gamma (-\nu -\frac d2)\Gamma (-\nu -d)
\end{equation}
so 
\begin{equation}
a_\nu ^1=b_{-\nu }^1=(2\pi )^{\frac{3d}2}\Gamma (-\nu )\Gamma (-\nu -\frac
d2)\Gamma (-\nu -d)
\end{equation}
according to the behaviour of the solutions $J_\nu ^{\left[ 3,j\right] }$.
To determine the other coefficients we take $x\neq 0$ on the boundary of $%
\Omega $. So if $x=x_1c_1+$ $x_2c_2$ then the integral representation of $%
K_\nu ^{\left[ 3\right] }$ takes the explicit form

\begin{eqnarray}
K_\nu ^{\left[ 3\right] }(x_1c_1+x_2c_2) &=&\stackunder{\Omega _2}{\dint }%
\stackunder{E}{\diint }\exp -\left[ \frac{2\det (z)+2ttr(z)-\left\| \xi
\right\| ^2}{2\det (z)\left[ t-\left( z^{-1}\xi ,\xi \right) \right] }%
\right] \;\times  \label{aut1} \\
&&\ \times \exp -\left( x,z\right) \;\det (z)^{\nu -d-1}\left[ t-\left(
z^{-1}\xi ,\xi \right) \right] ^{\nu -d-1}\;dzd\xi dt  \nonumber
\end{eqnarray}
where $E=\left\{ (t,\xi )\in \TeXButton{R}{{\R}}\times A_{1/2}\;/\;t>\left(
z^{-1}\xi ,\xi \right) \right\} $ and $y=z+\xi +tc_3$. We change $t$ by $%
t+\left( z^{-1}\xi ,\xi \right) $. The integral over $E$ becomes 
\begin{eqnarray*}
I &=&\stackrel{+\infty }{\stackunder{0}{\int }}\stackunder{A_{1/2}}{\dint }%
e^{-1/t}\exp -\left[ \frac{tr(z)\left( z^{-1}\xi ,\xi \right) -\frac
12\left\| \xi \right\| ^2}{t\det (z)}\right] \;t^{\nu -d-1}\;d\xi dt \\
\ &=&\stackrel{+\infty }{\stackunder{0}{\int }}\stackunder{A_{1/2}}{\dint }%
e^{-1/t}e^{-\left( B.\xi ,\xi \right) }\;t^{\nu -d-1}\;d\xi dt
\end{eqnarray*}
where the operator $B=\rho (v)$ with $v=\dfrac{tr(z)}{2t\det (z)}%
\;z^{-1}-\dfrac 1{2t\det (z)}\;e_0$. Note that $v\in \Omega _2$ because $%
tr(v)=\dfrac{tr(z)^2-2\det (z)}{2t\left( \det (z)\right) ^2}>0$ and $\det
(v)=\dfrac 1{4t^2\left( \det (z)\right) ^2}>0$. Also $\det B=\det
(v)^d=2^{-2d}t^{-2d}\det (z)^{-2d}.$ But 
\[
\stackunder{A_{1/2}}{\dint }e^{-\left( B.\xi ,\xi \right) }d\xi =\pi
^d\left( \det B\right) ^{-1/2} 
\]
so 
\begin{equation}
I=\left( 2\pi \right) ^d\det (z)^d\stackrel{+\infty }{\stackunder{0}{\int }}%
e^{-1/t}t^{\nu -1}dt=\left( 2\pi \right) ^d\Gamma (-\nu )\det (z)^d
\label{int}
\end{equation}
Now 
\begin{eqnarray}
K_\nu ^{\left[ 3\right] }(x_1c_1+x_2c_2) &=&\left( 2\pi \right) ^d\Gamma
(-\nu )\stackunder{\Omega _2}{\dint }e^{-\dfrac{tr(z)}{\det (z)}}\exp
-\left( x,z\right) \;\det (z)^{\nu -1}\;dz  \nonumber \\
\ &=&\left( 2\pi \right) ^d\Gamma (-\nu )\stackunder{\Omega _2}{\dint }%
e^{-tr(z^{-1})}\exp -\left( x,z\right) \;\det (z)^{\nu +d/2-d/2-1}\;dz 
\nonumber \\
\ &=&\left( 2\pi \right) ^d\Gamma (-\nu )K_{\nu +d/2}^{\left[ 2\right]
}(x_1c_1+x_2c_2)  \label{k3}
\end{eqnarray}

\begin{theorem}
We have (according to $\left( \ref{k3}\right) $,$\left( \ref{red1}\right) $,$%
\left( \ref{red2}\right) $,$\left( \ref{red3}\right) $,$\left( \ref{red4}%
\right) $ and the behaviour of J-solutions): 
\[
K_\nu ^{\left[ 3\right] }(x)=\stackunder{j=1}{\stackrel{4}{\dsum }}a_\nu
^jJ_\nu ^{\left[ 3,j\right] }+b_\nu ^jt_3^{-\nu }J_{-\nu }^{\left[
3,j\right] } 
\]
with 
\[
a_\nu ^1=b_{-\nu }^1=(2\pi )^{\frac{3d}2}\Gamma (-\nu )\Gamma (-\nu -\frac
d2)\Gamma (-\nu -d) 
\]
\[
a_\nu ^2=b_{-\nu }^2=(2\pi )^{\frac{3d}2}\Gamma (-\nu )\Gamma (-\nu -\frac
d2)\Gamma (\nu +d) 
\]
\[
a_\nu ^3=b_{-\nu }^3=(2\pi )^{\frac{3d}2}\Gamma (-\nu )\Gamma (\nu +\frac
d2)\Gamma (\nu ) 
\]
\[
a_\nu ^4=b_{-\nu }^4=(2\pi )^{\frac{3d}2}\Gamma (-\nu )\Gamma (\nu +\frac
d2)\Gamma (-\nu ) 
\]
\end{theorem}

\end{document}